\documentclass{amsart}
\usepackage{graphicx}
\usepackage{amssymb}
\usepackage{amsmath}
\usepackage{rotating}
\usepackage{setspace}
\usepackage{adjustbox}
\newtheorem{teor}{Theorem}
\newtheorem{ex}{Example}
\newtheorem{lem}[teor]{Lemma}
\newcommand{\be}{\begin{eqnarray}}
\newcommand{\ee}{\end{eqnarray}}
\newcommand{\bc}{\begin{center}}
\newcommand{\ec}{\end{center}}
\newcommand{\ds}{\displaystyle}
\newenvironment{demo}
{{\bf Proof.} }

\begin{document}

\title[Well-posedness and numerical approximation of tempered fractional terminal value problems]{Well-posedness and numerical approximation of tempered fractional terminal value problems}%

\author{M.~L.~Morgado}

\address{Centre of Mathematics, Pole CMAT-UTAD and Department of Mathematics, University of Tr\'as-os-Montes e Alto Douro, UTAD,
       Quinta de Prados 5001-801, Vila Real, Portugal}
       \email{luisam@utad.pt}
       \author{M.~Rebelo}
\address{Department of Mathematics and Centro de Matem\'atica e Aplica\c c\~oes, Universidade NOVA de Lisboa, Quinta da Torre, 2829-516, Caparica, Portugal}
 \email{msjr@fct.unl.pt}

\thanks{
Please always cite to this paper as: submitted to Fract. Calc. Appl. Anal., https://www.degruyter.com/view/j/fca and check there for further publication details.}%


\begin{abstract}
For a class of tempered fractional terminal value problems of the Caputo type, we study the existence and uniqueness of the solution, analyse the continuous dependence on the given data and using a shooting method, we present and discuss three numerical schemes for the numerical approximation of such problems. Some numerical examples are considered in order to illustrate the theoretical results and evidence the efficiency of the numerical methods.
\end{abstract}
\maketitle
\section{Introduction}
In this work we  analyse a class of tempered terminal value problems for tempered fractional ordinary differential equations of order $\alpha$, with $0<\alpha<1$:
\be && \label{eq11}
{}_0\mathbb{D}_t^{\alpha,\lambda}
\left(y(t)\right)=f(t,y(t)),\quad t\in [0,a],\\
 && \label{eq12} e^{\lambda a}y(a)=y_a,
\ee
where $f$ is a suitably behaved function and $\ds \mathbb{D}_0^{\alpha,\lambda}\left(y(t)\right)$  denotes the left-sided Caputo tempered fractional derivative of order $\alpha>0$, where the tempered parameter $\lambda$ is nonnegative.

The left-sided Caputo tempered fractional derivative can be given trough the definition of the Caputo derivative (see \cite{artigoArxiv} for example). In the particular case where  $0<\alpha<1$ it reads:
\be \label{defDerivada}
{}_0\mathbb{D}_t^{\alpha,\lambda}\left(y(t)\right)=
e^{-\lambda t}
{}_0\mathcal{D}_t^{C,\alpha}\left(e^{\lambda t} y(t)\right)
=
\frac{e^{-\lambda t} }{\Gamma(1-\alpha)}\int_{0}^{t}\frac{1}{(t-s)^{\alpha}}\frac{d\left(e^{\lambda s}y(s)\right)}{d s} ds, \ee where $\ds _0\mathcal{D}_t^{C,\alpha}$ denotes the Caputo fractional derivative (see \cite{Diethelm_2010}).

Note that if $\lambda=0$ then the Caputo tempered fractional derivative reduces to the Caputo fractional derivative, and therefore, Caputo derivatives can be regarded as a particular case of Caputo tempered derivatives.

Fractional differential equations of Caputo-type have been investigated extensively in the last decades and many significant contributions were provided by researchers of several areas as mathematics, physics and engineering, making Fractional Calculus as one of the most hot and current research topics. Recently, some attention has been devoted to tempered fractional differential equations, because the later ones revealed to model more realistically some phenomena (see \cite{Liemert} and the references therein for details). Even so, the literature is not so vast for this type of equations, as it is for fractional differential equations in the Caputo sense. As it happens with non-tempered fractional  differential equations, the analytical solution is usually impossible to obtain and in the cases where it can be determined, its representation in terms of a series makes it difficult to handle. Therefore, the development of numerical methods for this type of fractional differential equations is also crucial. With this respect, some approaches have already been reported. In \cite{Baeumer}, the authors propose a finite difference formula for tempered fractional derivatives and introduce a temporal and spatial second-order Crank-Nicolson method for the space-fractional diffusion equation. In \cite{artigoArxivAA} and \cite{artigoArxiv} a Jacobi-predictor-corrector algorithm is presented for tempered ordinary initial value problems. The authors in \cite{Marom} present a finite difference scheme to solve fractional partial differential models in finance. In \cite{Zhao}, spectral methods are derived for the tempered advection and diffusion problems.

To the best of our knowledge, tempered terminal value problems have never been investigated. Therefore, after analysing the well-posedness of such problems, we consider a simple approach for the numerical approximation of the solution, which is based on the relationship between tempered and non-tempered Caputo derivatives.\\
The paper is organized in the following way: in the next section we establish sufficient conditions for the existence and uniqueness of the solution of problems of the type (\ref{eq11})-(\ref{eq12}). Then we investigate the continuous dependence of the solution on the given data. Section \ref{Methods} is devoted to the derivation of numerical schemes and finally in section \ref{NumEx} we present and discuss several numerical examples. The paper ends with some conclusions and plans for further investigation.

\section{Existence and uniqueness of the solution } \label{sec2}

From \cite{artigoArxiv} we have the following two results for initial value problems.

\begin{lem}\label{lema1}\cite{artigoArxiv} If the function $f(t,u)$ is continuous, then the initial value problem
\be \label{IVP}
\left\{
\begin{array}{l}
{}_0\mathbb{D}_t^{\alpha,\lambda}
\left(y(t)\right)=f(t,y(t)),\quad t\in [0,a],\\
\left(e^{-\lambda t}y(t)\right)|_{t=0}=c_k,\, k=0,1,\ldots,n-1,
\end{array}
\right.
\ee
 is equivalent to the nonlinear Volterra integral equation of the second kind
\be\label{eqint2}
y(t)=\sum_{k=0}^{n-1}c_k\frac{e^{-\lambda t}t^k}{\Gamma(k+1)}+\frac{1}{\Gamma(\alpha)}\int_0^te^{-\lambda(t-s)}(t-s)^{\alpha-1}f(s,y(s))ds,
\ee
where $\ds n-1<\alpha<n$ and $\ds \alpha \ge 0$.\\
In the particular case where $0< \alpha <1$, we have
\[y(t)=c_0e^{-\lambda t}+\frac{1}{\Gamma(\alpha)}\int_0^te^{-\lambda(t-s)}(t-s)^{\alpha-1}f(s,y(s))ds.\]
\end{lem}
\begin{teor}\label{Teorema IVP}
Let  $\ds n-1<\alpha<n$, $n\in \mathbb{N}^+$, $\lambda \geq 0$ and $a\in \mathbb{R}$ such that $a>0$, then the fractional differential equation (\ref{IVP}) has a unique solution $u(t)\in C^n([0,a])$.
\end{teor}

Next, we extend these results to terminal boundary problems.
In what follows  the Caputo tempered fractional derivative of order $\alpha$ will be simply denoted by $\ds \mathbb{D}^{\alpha,\lambda}
\left(y(t)\right)$.

The Caputo tempered fractional derivative of order $\alpha$, with $\alpha \in (0,1)$ satisfies
\be&&\label{prop1}
\mathbb{I}^{\alpha,\lambda}\left(\mathbb{D}^{\alpha,\lambda}(u(t))\right)= u(t)-e^{\lambda t}u(0),
\\
&&\label{prop2}\mathbb{D}^{\alpha,\lambda}\left(\mathbb{I}^{\alpha,\lambda}(u(t))\right)=u(t),
\ee
where  $\ds \mathbb{I}^{\alpha,\lambda}$ is the Riemann-Liouville tempered fractional integral given by
\be\label{integral}
\mathbb{I}^{\alpha,\lambda}(u(t))=e^{-\lambda t}\mathbb{I}^{\alpha}(e^{\lambda t }u(t))=\frac{1}{\Gamma(\alpha)}\int_0^t e^{-\lambda(t-s)}(t-s)^{\alpha-1}u(s)ds,
\ee
and $\mathbb{I}^{\alpha}$ denotes the Riemann-Liouville fractional integral.

If $y$ satisfies the fractional differential equation (\ref{eq11}), then applying the Riemann-Liouville fractional integral at the both sides of equation and taking property (\ref{prop1}) into account, we conclude that the solution $y$ satisfies the following integral equation
\be\nonumber e^{\lambda t} y(t)&=&y(0)+\frac{1}{\Gamma(\alpha)}\int_0^{t}e^{\lambda s}
(t-s)^{\alpha-1}f(s,y(s))ds\\
\nonumber &=&y(0)+\frac{ 1}{\Gamma(\alpha)}\int_0^{a}e^{\lambda s}
(a-s)^{\alpha-1}f(s,y(s))ds\\
\nonumber
&+&\frac{1}{\Gamma(\alpha)}\left(-\int_0^{a}e^{\lambda s}
(a-s)^{\alpha-1}f(s,y(s))ds+\int_0^{t}e^{\lambda s}
(t-s)^{\alpha-1}f(s,y(s))ds\right)\\
\nonumber &=&y_a-\frac{1}{\Gamma(\alpha)}\left(\int_0^{a}e^{\lambda s}
(a-s)^{\alpha-1}f(s,y(s))ds-\int_0^{t}e^{\lambda s}
(t-s)^{\alpha-1}f(s,y(s))ds\right).\\
&&\label{VIeq}
\ee

Therefore,  if $y$ is a vsolution of the fractional boundary value problem (FBVP) (\ref{eq11})-(\ref{eq12}) then $y$ is a solution of the integral equation
\be\nonumber y(t)=y_a e^{-\lambda t }-\frac{e^{-\lambda t }}{\Gamma(\alpha)}\left(\int_0^{a}e^{\lambda s}
(a-s)^{\alpha-1}f(s,y(s))ds-\int_0^{t}e^{\lambda s}
(t-s)^{\alpha-1}f(s,y(s))ds\right)\\
&& \label{VIE}.\ee

Next, we  establish sufficient conditions
for the existence and uniqueness of solutions of the FBVP (\ref{eq11})-(\ref{eq12}). The proof will be based on the Banach's fixed point theorem.
We just establish the existence and uniqueness of the solution  on the interval $[0,a]$, since the existence and uniqueness for $t>a$ is inherited from the corresponding initial value problem theory (see\cite{artigoArxiv} for details).

Define the set $\ds \Omega_{\gamma}=\{y\in C([0,a]) : \, \|y-y_ae^{-\lambda t}\|_{[0,a]}\leq \gamma\}$, where the norm $\ds \|\cdot \|_{[0,a]}$  is defined by $\ds \|\cdot \|_{[0,a]}=\max_{t\in[0,a]}|g(t)|$, for all $g\in C([0,a])$ and
\be \label{defigamma}  \gamma=\frac{2a^{\alpha}\|f\|_{[0,a]}e^{\lambda a}}{\Gamma(1+\alpha)
}.\ee The set $\Omega_{\gamma}$ is a closed subset of
the Banach space of all continuous functions on $[0,a]$, equipped with the
norm $\ds \| \quad \|_{[0,a]}$, and since the function $\ds y(t)=y_ae^{-\lambda t} \in \Omega_{\gamma}$,  it is nonempty. On $\Omega_{\gamma}$, let us define
the operator
\be\nonumber (\mathcal{A} y)(t)&=& y_a e^{-\lambda t }-\frac{1}{\Gamma(\alpha)}\int_0^{a}e^{-\lambda (t-s)}
(a-s)^{\alpha-1}f(s,y(s))ds\\
\label{defoperador} &+&\int_0^{t}e^{-\lambda (t-s)}
(t-s)^{\alpha-1}f(s,y(s))ds.\ee
Using this operator, the integral equation can be rewritten as $\ds y=\mathcal{A} y$, and if the operator $\mathcal{A}$ has a unique fixed point on $U$ then the FBVP (\ref{eq11})-(\ref{eq12}) has a unique continuous solution.  Using the Banach's fixed point theorem, under some assumptions on $\ds f$, we prove the existence and uniqueness result on the next theorem.

\begin{teor}\label{Teorema1}
Let   $\ds D=[0,a]\times [y_ae^{-\lambda a }-\gamma,y_ae^{-\lambda a }+\gamma]$, with $\gamma$ given by (\ref{defigamma}), and assume that  the function $f:D\rightarrow \mathbb{R}$ is continuous for all $t\in [0,a].$
 We further assume that the function $f$ fulfills a Lipschitz condition with respect to the second variable, meaning that there exists $L>0$ such that it holds
 \be\label{LC}|f(t,y)-f(t,z)|\le L|y-z|, ~~\mbox{for}~~ y,\,z\in \Omega_{\gamma}.  \ee
 If the Lipschitz constant $L$ is such that $\ds L<\frac{\Gamma(\alpha+1)}{2a^{\alpha}e^{\lambda a}}$, then $\mathcal{A}$ maps $\Omega_{\gamma}$ into itself and it is a contraction:
\be\label{contract}
\|\mathcal{A}y-\mathcal{A}z\|_{[0,a]}\le \|y-z\|_{[0,a]}
\quad\text{for}\quad y,z\in \Omega_{\gamma}.
\ee
Hence equation (\ref{VIE}) has a unique solution $y^*\in\Omega_{\gamma}$, which is  the unique fixed point of $\mathcal{A}.$
\end{teor}

\begin{demo}

Let $\ds y\in \Omega_{\gamma}$. First, we show that $\mathcal{A}y\in \Omega_{\gamma}$. \\
From the definition of $\mathcal{A}$  we have
\be \nonumber \left|\mathcal{A}y-y_ae^{-\lambda t}\right|&=&\frac{e^{-\lambda t}}{\Gamma(\alpha)}\left|
\int_0^{a}e^{\lambda s}
(a-s)^{\alpha-1}f(s,y(s))ds+\int_0^{t}e^{\lambda s}
(t-s)^{\alpha-1}f(s,y(s))ds \right|\\
\nonumber &\leq &
\frac{\|f\|_{[0,a]}e^{\lambda T}}{\alpha\Gamma(\alpha)
}(a^{\alpha}+t^{\alpha}).\ee
Then
\be \left\|\mathcal{A}y-y_ae^{-\lambda t}\right\|_{[0,a]}\leq
\frac{2 a^{\alpha}\|f\|_{[0,a]}e^{\lambda a}}{\Gamma(\alpha+1)
}=\gamma,\ee which implies that $\ds (\mathcal{A}y)\in \Omega_{\gamma}.$

Now we   prove that $\mathcal{A}$ is a contraction on $\Omega_{\gamma}$, with $\gamma$ defined by (\ref{defigamma}).
For $y,z\in \Omega_{\gamma}$, $t\in[0,a]$, we have
\be\nonumber
|(\mathcal{A}y)(t)-(\mathcal{A}z)(t)|&
\le& \frac{1}{\Gamma(\alpha)}\left(\int_{0}^{a}e^{-\lambda(t-s)}(a-s)^{\alpha-1}|f(s,y(s))-f(s,z(s))|ds\right.\\
\nonumber &+&\left.
\int_{0}^{t}e^{-\lambda(t-s)}(t-s)^{\alpha-1}|f(s,y(s))-f(s,z(s))|ds\right) \\
\nonumber
&\le&\frac{2La^{\alpha}e^{\lambda a}}{\alpha\Gamma(\alpha)}\|y-z\|_{[0,a]}=\frac{2La^{\alpha}e^{\lambda a}}{\Gamma(1+\alpha)}\|y-z\|_{[0,a]}
< \|y-z\|_{[0,a]}.\ee
 Then, the operator  $\mathcal{A}$ is a contraction on $\Omega_{\gamma}$.  Finally, by the Banach fixed point principle the proof of the theorem is complete.
\end{demo}

\vspace*{0.5cm}

If the assumptions of Theorem \ref{Teorema1} are satisfied, then the FBVP (\ref{eq11})-(\ref{eq12}) has a unique continuous solution, $y(t)$, on the interval $[0,a]$ and, in particular, a unique value for $y(0)$  exists. Therefore, there is an exact correspondence between tempered fractional
boundary value problems and tempered fractional initial value problems.

\section{Continuous dependence of the solution on the data}
\
\noindent

In order to analyse the continuous dependence of the solution on the given data we assume that problem
\begin{eqnarray}  \nonumber
\mathbb{D}_t^{\alpha,\lambda}
\left(y(t)\right)&=&f(t,y(t)),\quad t\in [0,a],\\
\nonumber e^{\lambda a}y(a)&=&y_a,
\end{eqnarray}
which is equivalent to
\be\nonumber &&
y(t)=y_a e^{-\lambda t }-\frac{e^{-\lambda t }}{\Gamma(\alpha)}\int_0^{a}e^{\lambda s}
(a-s)^{\alpha-1}f(s,y(s))ds+\frac{e^{-\lambda t }}{\Gamma(\alpha)}\int_0^{t}e^{\lambda s}
(t-s)^{\alpha-1}f(s,y(s))ds,\\
\label{equiv1}\ee\emph{•}
may suffer perturbations on the parameters $y_a$, $\alpha$, $\lambda$ and on the right-hand side function $f$, and therefore we will consider the following perturbed problems:
\begin{eqnarray} \label{eq11a}
\mathbb{D}_t^{\alpha,\lambda}
\left(z(t)\right)&=&f(t,z(t)),\quad t\in [0,a],\\
 \label{eq12a} e^{\lambda a}z(a)&=&z_a,
\end{eqnarray}

\begin{eqnarray} \label{eq11b}
\mathbb{D}_t^{\alpha-\delta,\lambda}
\left(z(t)\right)&=&f(t,z(t)),\quad t\in [0,a],\\
 \label{eq12b} e^{\lambda a}z(a)&=&y_a,
\end{eqnarray}

\begin{eqnarray} \label{eq11c}
\mathbb{D}_t^{\alpha,\lambda-\delta}
\left(z(t)\right)&=&f(t,z(t)),\quad t\in [0,a],\\
 \label{eq12c} e^{\left(\lambda -\delta\right)a}z(a)&=&y_a,
\end{eqnarray}

\begin{eqnarray} \label{eq11d}
\mathbb{D}_t^{\alpha,\lambda}
\left(z(t)\right)&=&\tilde{f}(t,z(t)),\quad t\in [0,a],\\
 \label{eq12d} e^{\lambda a}z(a)&=&y_a,
\end{eqnarray}
$\delta>0$, and are equivalent to the following integral equations:

\be \nonumber &&
z(t)=z_a e^{-\lambda t }-\frac{e^{-\lambda t }}{\Gamma(\alpha)}\int_0^{a}e^{\lambda s}
(a-s)^{\alpha-1}f(s,z(s))ds+\\
\label{equiv2} &&\quad \quad \quad +\frac{e^{-\lambda t }}{\Gamma(\alpha)}\int_0^{t}e^{\lambda s}
(t-s)^{\alpha-1}f(s,z(s))ds,\\
\nonumber && \\
\nonumber &&
z(t)=y_a e^{-\lambda t }-\frac{e^{-\lambda t }}{\Gamma(\alpha-\delta)}\int_0^{a}e^{\lambda s}
(a-s)^{\alpha -\delta-1}f(s,z(s))ds+\\
\label{equiv3}&&\quad \quad \quad +\frac{e^{-\lambda t }}{\Gamma(\alpha-\delta)}\int_0^{t}e^{\lambda s}
(t-s)^{\alpha -\delta-1}f(s,z(s))ds,\\
\nonumber &&\\
\nonumber &&z(t)=y_a e^{-\left(\lambda -\delta\right)t }-\frac{e^{-\left(\lambda-\delta\right) t }}{\Gamma(\alpha)}\int_0^{a}e^{\left(\lambda-\delta\right) s}
(a-s)^{\alpha-1}f(s,z(s))ds+
\\
\label{equiv4}&&\quad \quad \quad +\frac{e^{-\left(\lambda-\delta\right) t }}{\Gamma(\alpha)} \int_0^{t}e^{\left(\lambda -\delta\right) s}
(t-s)^{\alpha-1}f(s,z(s))ds,\\
\nonumber &&\\
\nonumber &&y(t)=y_a e^{-\lambda t }-\frac{e^{-\lambda t }}{\Gamma(\alpha)}\int_0^{a}e^{\lambda s}
(a-s)^{\alpha-1}\tilde{f}(s,z(s))ds+\\
\label{equiv5}&&\quad \quad \quad+\frac{e^{-\lambda t }}{\Gamma(\alpha)}-\int_0^{t}e^{\lambda s}
(t-s)^{\alpha-1}\tilde{f}(s,z(s))ds,
\ee
respectively.
\begin{teor}\label{Varua}
Let $y$ and $z$ be the unique solutions of problems (\ref{eq11})-(\ref{eq12}) and (\ref{eq11a})-(\ref{eq12a}), respectively. Then
\[\left\|y-z\right\| \le \frac{1}{\beta}\left|y_a-z_a\right|,\]
where $\beta=1-\frac{2Le^{\lambda a}a^{\alpha}}{\Gamma(\alpha+1)}$.
\end{teor}
\begin{demo}
Taking (\ref{equiv1}) and (\ref{equiv2}) into account, for any $t \in [0,a]$, we have
\begin{eqnarray}\nonumber \left|y(t)-z(t)\right|  && \le  e^{-\lambda t }\left|\left(y_a-z_a\right)-\frac{1}{\Gamma(\alpha)}\int_0^{a}e^{\lambda s}(a-s)^{\alpha-1}\left(f(s,y(s))-f(s,z(s))\right)ds\right|+\\
\nonumber &&\frac{e^{-\lambda t }}{\Gamma(\alpha)}\int_0^{t}e^{\lambda s}(t-s)^{\alpha-1}\left|f(s,y(s))-f(s,z(s))\right|ds\\
\nonumber && \le  \left|\left(y_a-z_a\right)-\frac{1}{\Gamma(\alpha)}\int_0^{a}e^{\lambda s}
(a-s)^{\alpha-1}\left(f(s,y(s))-f(s,z(s))\right)ds\right|+\\
\nonumber && \frac{Le^{\lambda a}}{\Gamma(\alpha)}\int_0^{t}
(t-s)^{\alpha-1}\left|y(s)-z(s)\right|ds\\
\nonumber && \le \left|y_a-z_a\right|+\frac{Le^{\lambda a}}{\Gamma(\alpha)}\int_0^{a}
(a-s)^{\alpha-1}\left|y(s)-z(s)\right|ds\\
\nonumber && +\frac{Le^{\lambda a}}{\Gamma(\alpha)}\int_0^{t}
(t-s)^{\alpha-1}\left|y(s)-z(s)\right|ds.
\end{eqnarray}
Hence
\begin{eqnarray}
\nonumber \left\|y-z\right\|&& \le \left|y_a-z_a\right|+\frac{Le^{\lambda a}}{\Gamma(\alpha)}\left\|y-z\right\|\left(\int_0^{a}(a-s)^{\alpha-1}ds+\int_0^{t}(t-s)^{\alpha-1}ds\right)\\
\nonumber && =\left|y_a-z_a\right|+\frac{Le^{\lambda a}}{\Gamma(\alpha)}\
\left\|y-z\right\|\left(\frac{a^{\alpha}}{\alpha}+\frac{t^{\alpha}}{\alpha}\right)\\
\nonumber && \le \left|y_a-z_a\right|+\frac{2Le^{\lambda a}a^{\alpha}}{\Gamma(\alpha+1)}\
\left\|y-z\right\|.
\end{eqnarray}
According to the upper bound on the Lipschitz condition $L$, established in Theorem 3, we have
\begin{equation}\label{beta}\beta=1-\frac{2Le^{\lambda a}a^{\alpha}}{\Gamma(\alpha+1)}>0,\end{equation}
and therefore, we conclude that
\[\left\|y-z\right\| \le \frac{1}{\beta}\left|y_a-z_a\right|,\]
and the Theorem is proved.
\end{demo}
\begin{teor}\label{Varalpha}
Let $y$ and $z$ be the unique solutions of problems (\ref{eq11})-(\ref{eq12}) and (\ref{eq11b})-(\ref{eq12b}), respectively, where in the later we assume that $\delta$ is such that \\ $0< \alpha -\delta <1$. Then
\[\left\|y-z\right\| =\mathcal{O}\left(\delta \right).\]
\end{teor}
\begin{demo}
Taking (\ref{equiv1}) and (\ref{equiv3}) into account, for any $t \in [0,a]$, we have
\begin{eqnarray}\nonumber \left|y(t)-z(t)\right|  &=&  \left|-\frac{e^{-\lambda t }}{\Gamma(\alpha)}\int_0^{a}e^{\lambda s}(a-s)^{\alpha-1}f(s,y(s))ds +\right.\\
\nonumber &+&\left.\frac{e^{-\lambda t }}{\Gamma(\alpha -\delta)}\int_0^{a}e^{\lambda s}(a-s)^{\alpha-\delta-1}f(s,z(s))ds \right.\\
\nonumber && \left.+ \frac{e^{-\lambda t }}{\Gamma(\alpha)}\int_0^{t}e^{\lambda s}(t-s)^{\alpha-1}f(s,y(s))ds\right.\\
\nonumber  &-&\left.\frac{e^{-\lambda t }}{\Gamma(\alpha -\delta)}\int_0^{t}e^{\lambda s}(t-s)^{\alpha-\delta-1}f(s,z(s))ds\right|
\end{eqnarray}
Since $e^{-\lambda t} \le 1$, for all $t\ge 0$ and $\lambda >0$ and $e^{\lambda s} \le e^{\lambda a}$, for all $0 \le s \let \le a$, then
\begin{eqnarray}\nonumber \left|y(t)-z(t)\right|  & \le& e^{\lambda a}\left(\int_0^{a}\left|\frac{(a-s)^{\alpha-1}}{\Gamma(\alpha)}f(s,y(s))-\frac{(a-s)^{\alpha-\delta-1}}{\Gamma(\alpha-\delta)}f(s,z(s))\right|~ds+\right.\\
\label{eqaux} &+& \left.\int_0^{t}\left|\frac{(t-s)^{\alpha-1}}{\Gamma(\alpha)}f(s,y(s))-\frac{(t-s)^{\alpha-\delta-1}}{\Gamma(\alpha-\delta)}f(s,z(s))\right|~ds\right).
\end{eqnarray}
Let us look at the first integral on the right-hand side of (\ref{eqaux}) first:
\begin{eqnarray}
\nonumber I_1&=& \int_0^{a}\left|\frac{(a-s)^{\alpha-1}}{\Gamma(\alpha)}f(s,y(s))-\frac{(a-s)^{\alpha-\delta-1}}{\Gamma(\alpha-\delta)}f(s,z(s))\right|~ds\\
\nonumber &=& \int_0^{a}\left|\frac{(a-s)^{\alpha-1}}{\Gamma(\alpha)}f(s,y(s))-\frac{(a-s)^{\alpha-1}}{\Gamma(\alpha)}f(s,z(s))\right.\\
\nonumber &&\quad \quad  +\left.\frac{(a-s)^{\alpha-1}}{\Gamma(\alpha)}f(s,z(s))-\frac{(a-s)^{\alpha-\delta-1}}{\Gamma(\alpha-\delta)}f(s,z(s))\right|~ds\\
\nonumber &\le & \frac{La^{\alpha}}{\Gamma(\alpha +1)}\left\|y-z\right\|+\int_0^{a}\left|\frac{(a-s)^{\alpha-1}}{\Gamma(\alpha)}-\frac{(a-s)^{\alpha-\delta-1}}{\Gamma(\alpha-\delta)}\right|\left|f(s,z(s))\right|ds.
\end{eqnarray}
Considering the function $\displaystyle G(x)=\frac{(a-s)^{x}}{\Gamma(x+1)}$, using the mean value Theorem, we easily conclude that $ \displaystyle\left|\frac{(a-s)^{\alpha-1}}{\Gamma(\alpha)}-\frac{(a-s)^{\alpha-\delta-1}}{\Gamma(\alpha-\delta)}\right| \le C \delta$, where $\displaystyle C=\max_{x\in[\alpha-\delta-1,\alpha-\delta]}G'(x)$, and therefore
\[I_1 \le \frac{La^{\alpha}}{\Gamma(\alpha +1)}\left\|y-z\right\|+Ca\left\|f\right\|\delta.\]
Proceeding similarly with the second integral on the right-hand side of (\ref{eqaux}), we could conclude that
\be \nonumber I_2&=&\int_0^{t}\left|\frac{(t-s)^{\alpha-1}}{\Gamma(\alpha)}f(s,y(s))-\frac{(t-s)^{\alpha-\delta-1}}{\Gamma(\alpha-\delta)}f(s,z(s))\right|~ds\\
\nonumber&\le&
\frac{Lt^{\alpha}}{\Gamma(\alpha +1)}\left\|y-z\right\|+Ct\left\|f\right\| \delta \le \frac{La^{\alpha}}{\Gamma(\alpha +1)}\left\|y-z\right\|+Ca\left\|f\right\|\delta,\ee
and therefore
\[\left\|y-z\right\| \le \frac{2La^{\alpha}e^{\lambda a}}{\Gamma(\alpha +1)}\left\|y-z\right\|+2Ca\left\|f\right\|\delta,\]
or
\[\left\|y-z\right\| \le \frac{2Ca\left\|f\right\|\delta}{\beta},\]
where $\beta$ is given by (\ref{beta}). This completes the proof of the Theorem.
\end{demo}
\begin{teor}\label{Varlambda}
Let $y$ and $z$ be the unique solutions of problems (\ref{eq11})-(\ref{eq12}) and (\ref{eq11c})-(\ref{eq12c}), respectively. Then
\[\left\|y-z\right\| =\mathcal{O}\left(\delta \right).\]
\end{teor}
\begin{demo}
Taking (\ref{equiv1}) and (\ref{equiv4}) into account, for any $t \in [0,a]$, we have
\begin{eqnarray}
\nonumber \left|y(t)-z(t)\right| &\le& y_a\left|e^{-\lambda t}-e^{-\left(\lambda-\delta\right) t}\right|+\frac{1}{\Gamma(\alpha)}\int_0^{a}\left|e^{\lambda s}f(s,y(s))-e^{\left(\lambda-\delta\right) s}f(s,z(s))\right|(a-s)^{\alpha-1}~ds\\
\label{eqaux2} &+&\frac{1}{\Gamma(\alpha)}\int_0^{t}\left|e^{\lambda s}f(s,y(s))-e^{\left(\lambda-\delta\right) s}f(s,z(s))\right|(t-s)^{\alpha-1}~ds.
\end{eqnarray}
By using the mean value Theorem with the function $g(x)=e^{-xt}$, we conclude that
\[\left|e^{-\lambda t}-e^{-\left(\lambda-\delta\right) t}\right|\le   C_1\delta, \]
where $\displaystyle C_1=\max_{x\in [\lambda -\delta,\lambda]}\left|g'(x)\right|$.\\
Concerning the first integral on the right-hand side of (\ref{eqaux2}):
\begin{eqnarray}
\nonumber J_1 &=& \int_0^{a}\left|e^{\lambda s}f(s,y(s))-e^{\left(\lambda-\delta\right) s}f(s,z(s))\right|(a-s)^{\alpha-1}~ds\\
\nonumber &\le & \int_0^{a}\left|e^{\lambda s}f(s,y(s))-e^{\lambda s}f(s,z(s))\right|(a-s)^{\alpha-1}~ds+\\
\nonumber &+&\int_0^{a}\left|e^{\lambda s}f(s,z(s))-e^{\left(\lambda-\delta\right) s}f(s,z(s))\right|(a-s)^{\alpha-1}~ds\\
\nonumber & \le & \frac{e^{\lambda a}La^{\alpha}}{\alpha}\left\|y-z\right\|+C_2\left\|f\right\|\delta,
\end{eqnarray}
where by the mean value Theorem, $\displaystyle C_2=\max_{x \in [\lambda-\delta,\lambda]}\left|h'(x)\right|$, where $\displaystyle h(x)=e^{x s}$.\\
Proceeding analogously with the second integral in (\ref{eqaux2}), we conclude that
\[J_2=\frac{1}{\Gamma(\alpha)}\int_0^{t}\left|e^{\lambda s}f(s,y(s))-e^{\left(\lambda-\delta\right) s}f(s,z(s))\right|(t-s)^{\alpha-1}~ds \le \frac{e^{\lambda a}La^{\alpha}}{\alpha}\left\|y-z\right\|+C_2\left\|f\right\|\delta.\]
Hence
\[\left\|y-z\right\| \le C_1y_a\delta+\frac{2e^{\lambda a}La^{\alpha}}{\Gamma(\alpha+1)}\left\|y-z\right\| +2C_2\left\|f\right\|\delta,\]
or
\[\left\|y-z\right\| \le \frac{C_1+2C_2\left\|f\right\|}{\beta}\delta,\]
with  $\beta$ defined in (\ref{beta}). Thus, the Theorem is proved.
\end{demo}
\begin{teor}\label{Varf}
Let $y$ and $z$ be the unique solutions of problems (\ref{eq11})-(\ref{eq12}) and (\ref{eq11d})-(\ref{eq12d}), respectively. Then
\[\left\|y-z\right\| =\mathcal{O}\left(\left\|f-\tilde{f}\right\| \right).\]
\end{teor}
\begin{demo}
Taking (\ref{equiv1}) and (\ref{equiv5}) into account, for any $t \in [0,a]$, we have
\be
\nonumber \left|y(t)-z(t)\right| &\le& \frac{e^{\lambda a}}{\Gamma(\alpha)}\int_0^{a}(a-s)^{\alpha-1}\left|f(s,y(s))-\tilde{f}(s,z(s))\right|~ds+\\
\label{eqauxf}&+&\frac{e^{\lambda a}}{\Gamma(\alpha)}\int_0^{t}(t-s)^{\alpha-1}\left|f(s,y(s))-\tilde{f}(s,z(s))\right|~ds.
\ee
Since
\begin{eqnarray}
\nonumber && \int_0^{a}(a-s)^{\alpha-1}\left|f(s,y(s))-\tilde{f}(s,z(s))\right|~ds \le\\
&& \nonumber =\int_0^{a}(a-s)^{\alpha-1}\left|f(s,y(s))-f(s,z(s))\right|~ds+\int_0^{a}(a-s)^{\alpha-1}\left|f(s,z(s))-\tilde{f}(s,z(s))\right|~ds\\
\nonumber && \le \frac{La^{\alpha}}{\alpha}\left\|y-z\right\|+\frac{a^{\alpha}}{\alpha}\left\|f-\tilde{f}\right\|,
\end{eqnarray}
and, analogously
\be \nonumber \int_0^{t}(t-s)^{\alpha-1}\left|f(s,y(s))-\tilde{f}(s,z(s))\right|~ds &\le& \frac{Lt^{\alpha}}{\alpha}\left\|y-z\right\|+\frac{t^{\alpha}}{\alpha}\left\|f-\tilde{f}\right\|\\
\nonumber &\le& \frac{La^{\alpha}}{\alpha}\left\|y-z\right\|+\frac{a^{\alpha}}{\alpha}\left\|f-\tilde{f}\right\|,\ee
then
\[\left\|y-z\right\| \le \frac{2La^{\alpha}e^{\lambda a}}{\Gamma(\alpha+1)}\left\|y-z\right\| +2\frac{a^{\alpha}e^{\lambda a}}{\Gamma(\alpha+1)}\left\|f-\tilde{f}\right\|,\]
and the result of the Theorem follows.
\end{demo}

\section{Numerical method}\label{Methods}
\
\noindent
The results proved in \cite{Diethelm2012} can be applied to the integral equation (\ref{eqint2}) for $\ds \alpha \in (0,1)$ and $a=0$:

\be\label{eqint2V2}
y(t)=c_0 e^{-\lambda t}+\frac{1}{\Gamma(\alpha)}\int_0^te^{-\lambda(t-s)}(t-s)^{\alpha-1}f(s,u(s))ds.
\ee

Indeed, we can rewrite the integral equation (\ref{eqint2}), with $\alpha \in (0,1)$ and $a=0$,  as
\be\label{VIEIVP2}
u(t)&=&c_0+\frac{1}{\Gamma(\alpha)}\int_0^t(t-s)^{\alpha-1}g(s,u(s))ds,
\ee
where $\ds c_0=y(0)$, and
\begin{equation}\label{fs}u(s)=e^{\lambda s}y(s)~~\mbox{and}~~\ds g(s,u(s))=e^{\lambda s}f(s,e^{-\lambda s}u(s)).\end{equation}
From Theorem 3.1 of  \cite{Diethelm2012}  follows the following result:
\begin{teor}\label{teor1VIE}
If the function $\ds f$ is continuous on $\ds [0,b]\times \mathbb{R}$ and is Lipschitz continuous with Lipschitz
constant $L$ with respect to its second argument, then
\begin{itemize}
\item[a)]  the integral equation (\ref{eqint2}) to the integral equation  has a unique continuous solution on $[0,b]$, where $b$ can be arbitrary large;
\end{itemize}
and
\begin{itemize}
\item[b)] for every $c \in \mathbb{R}$, there
exists precisely one value of $c_0 \in \mathbb{R}$ for which the solution $y$ of (\ref{eqint2}) satisfies
$y(b) = ce^{-\lambda b}$.
\end{itemize}
\end{teor}
From Theorem \ref{teor1VIE} it follows that if $y$ and $u$ satisfy the integrals equations (\ref{eqint2V2}) and  (\ref{VIEIVP2}) with $\ds c_0=y_0$ and $c_0=u_0$, respectively, then $\ds y(t)=u(t)$ for all $\ds t\in[0,b]$  if and only if $y_0=u_0.$ Therefore taking the equivalence between initial value tempered fractional equations  and integral equations into account (c.f. Theorem \ref{Teorema IVP}) we obtain the following result for   tempered fractional equations of order $\alpha$, with $\alpha \in(0,1)$.
\begin{teor}\label{teorunicidadesol}
Let $\alpha \in (0,1)$ and assume that $\ds f:[0,b]\times [c,d]\rightarrow \mathbb{R}$ is continuous and satisfies a Lipchitz condition with respect to the second variable.\\
If $y_1$ and $y_2$  are are two solutions of the tempered differential equations
\be\label{TFDEqs}\mathbb{D}^{\alpha,\lambda}
\left(y_j(t)\right)=f(t,y_j(t)),\quad j=1,2,
\ee
subject to the initial conditions $\ds y_j(0)=y_{j0},~j=1,2$, respectively, where $\ds y_{10}\neq y_{20}$.\\
Then for all $t$ where both $y_1(t)$ and $y_2(t)$ exists we have $y_1(t)\neq y_2(t).$
\end{teor}

From Theorem \ref{teorunicidadesol} we can conclude that a solution of a tempered fractional differential equation of order $\alpha \in(0,1)$ is uniquely defined by a condition that can be specified at any point $\ds t\geq 0$.\\
On the other hand, Theorem \ref{teorunicidadesol} will be crucial to properly define the ideas of the numerical methods that we present next.

From Theorem \ref{teorunicidadesol} it follows that for the solution of (\ref{eq11}) that passes through the point
$(a, \exp(-\lambda a)y_a)$, we are able to find at most one point $(0, y_0)$ that also lies on the same solution
trajectory.  In order to obtain an approximation of $y(0)$ we propose   a shooting algorithm based on the
bisection method.
Let $y_1$ and $y_2$ be the solutions of (\ref{TFDEqs}) with initial values $y_{01}$ and $y_{02}$, such that $\ds y_1(a) <\exp(-\lambda a)y_a<y_2(a) $, the bisection  method  provide successive approximations for $y(0)$
until the distance between the two last approximations does not exceed a
given tolerance $\epsilon$.\\
 To evaluate the value of $y(a)$  we need a numerical method to solve the initial value problems
\be && \label{eqIVP1}
\mathbb{D}^{\alpha,\lambda}
\left(y(t)\right)=f(t,y(t)),\quad t\in (0,a],\\
 && \label{eqIVP2} y(0)=y_0.
\ee
 This will be straightforward if we take relationship (\ref{defDerivada}) into account. In fact, defining the functions $u$
and $g$ as in (\ref{fs}), we can use any available solver for (non-tempered) Caputo-type initial value problems to determine the solution $u$ of
\be && \label{eqIVP1Cap}
{}_0\mathcal{D}_t^{C,\alpha}
\left(u(t)\right)=g(t,u(t)),\quad t\in (0,a],\\
 && \label{eqIVP2Cap} y(0)=y_0,
\ee
and then the solution of (\ref{eqIVP1})-(\ref{eqIVP2}) will be given by $\ds y(t)=e^{-\lambda t}u(t)$.

\section{Numerical results}\label{NumEx}
\
\noindent

In this section we present some numerical examples to illustrate the efficiency of the numerical algorithm.

\subsection{Approximating  the solution of the terminal value  problem (\ref{eq11})-(\ref{eq12})  }
\
\noindent

In this subsection we present 3 examples and the method that we apply for each one, depends on the nature of the differential equation and regularity of the solution.

The first example is a linear fractional differential equation with a smooth solution.
\begin{ex}\label{exemplo1}
\be\nonumber && \mathbb{D}^{\alpha,\lambda}
\left(y(t)\right)= e^{-\lambda  t} \left(\frac{3 \Gamma (3) t^{2-\alpha }}{4 \Gamma (3-\alpha )}+\frac{\Gamma (5) t^{4-\alpha }}{\Gamma (5-\alpha )}+c_{\alpha,\lambda}\left(t^4+\frac{3 t^2}{4}\right)\right)-c_{\alpha,\lambda} y(t),\quad t>0,\\
\nonumber &&
y(0.5)=\frac{e^{-\lambda/2}}{4},
\ee
\end{ex}
where $\ds c_{\alpha,\lambda}=\frac{\Gamma(\alpha+1)}{2^{1-\alpha}e^{\lambda/2}}$ and whose analytical solution is given by \\ $\ds y(t)=
\left(t^4+\frac{3 t^2}{4}\right) e^{ -\lambda  t}$ .


The second example is  a nonlinear fractional differential equation  with a smooth solution defined by:

\begin{ex}\label{exemplo2}
\be\nonumber && \mathbb{D}^{\alpha,\lambda}
\left(y(t)\right)= e^{-\lambda  t} \left(\frac{ \Gamma (3) t^{2-\alpha }}{ \Gamma (3-\alpha )}t^{2-\alpha}-3 t^4 \exp (-\lambda  t)\right)+3 y^2,\quad t>0,\\
\nonumber &&
y(0.5)=\frac{e^{-\lambda/2}}{4},
\ee
\end{ex}
 whose analytical solution is given by $\ds y(t)= e^{ -\lambda  t}t^2$ .


The third  example is  a linear fractional  differential equation  with a  solution whose second derivative has a singularity at $t=0$.

\begin{ex}\label{exemplo3}
\be\nonumber && \mathbb{D}^{\alpha,\lambda}
\left(y(t)\right)=
e^{-\lambda  t} \frac{ \Gamma (5/2) t^{2-\alpha }}{ \Gamma (5/2-\alpha )}t^{3/2-\alpha},\quad t>0,\\
\nonumber &&
y(0.5)=e^{-\lambda/2}\sqrt{\frac{1}{2^3}},
\ee
\end{ex}
 whose analytical solution is given by $\ds y(t)= e^{ -\lambda  t}t^{3/2}$ .

\vspace*{0.3cm}

In what follows we consider examples \ref{exemplo1},
\ref{exemplo2} and \ref{exemplo3}  with several values of $\alpha$  and $\lambda=2$.

 In order to compute $y(0)$  the bisection method was  used with  $\epsilon=10^{-10}$ and the approximate solution of each one of IVP was computed with the  three methods listed bellow.
\begin{itemize}
\item {\it{Method 1.}} Fractional backward difference based on quadrature (see,
for example, \cite{Diethelm_1997}).
\item {\it{Method 2.}} This numerical method can be seen as a generalization of the classical one-step Adams-Bashforth-Moulton scheme for first-order equations (cf. \cite{Diethelm_2004} ) and is appropriate to obtain a numerical solution  of the non-linear problems.

\item {\it{Method 3.}} In this method we consider an integral formulation of the initial value problem (\ref{eqIVP1Cap})-(\ref{eqIVP2Cap}) and a nonpolynomial approximation of the solution (cf.\cite{MorgadoFordRebelo2013} ). This method is appropriate to approximate the solution of  problems whose solution is not smooth.

\end{itemize}

We denote the absolute errors by $\ds e^h_{\epsilon}(t)=|y(t)-y^h_{\epsilon}(t)|$, where $y^h$ is an approximate solution of $y$ using the algorithm with  stepsize $\ds h=\frac{a}{N}$ and the value $y_0\sim y(0)$ was obtained by the bisection method with tolerance $\epsilon$.

 The absolute errors for
$\ds t=a=0.5$ and $\ds t = 1$ and the obtained values of $y(0)$ for the examples  \ref{exemplo1},
\ref{exemplo2} and \ref{exemplo3} are presented in Tables \ref{tab1}, \ref{tab3} and \ref{tab5exemp3}, respectively.\\
 In Tables \ref{tab2}, \ref{tab4} and \ref{tab6exemp3} the maximum of the absolute errors, \\
  $\|e^{a/N}\|_{\infty} =\ds \max_{0\leq i \leq N}e^{a/N}_{\epsilon}(t_i)$, and the experimental orders of convergence,  \\ $\ds p_N=\frac{\log(\|e^{a/N}\|_{\infty}/\|e^{a/2N}\|_{\infty})}{\log(2)}$, are listed.

In Table \ref{tab1} we observe that the absolute error for the point where the boundary condition is imposed, does not decrease as
the step-size goes smaller, although we are comparing very small quantities. On the other hand, for the approximate solution of Example \ref{exemplo2} the absolute  error at the boundary point decreases as the step-size $h$ decreases (cf. Table \ref{tab3}) and decreases with convergence order $1+\alpha$.

In Tables \ref{tab2} and  \ref{tab4} the experimental orders of convergence are listed, and we observe that the corresponding to examples \ref{exemplo1} and \ref{exemplo2} are approximately $2-\alpha$ and $1+\alpha$, respectively. The results  are in agreement with the theoretical result proved in  \cite{Diethelm_1997}, for Method 1, and  with the conjecture of Diethelm \textit{et al} \cite{Diethelm_2004}, for Method 2.

In Tables \ref{tab5exemp3} and \ref{tab6exemp3} we compare the results obtained  with the shooting method and $y_0$ given by Method 1 and Method 3 on the space $V_{h,1}^{1/2}$. We observe that the error at $t=0$ is smaller when $y_0$ is obtained by Method 3. However, both methods converge to zero with convergence order $1.5$, approximately. From  Figure \ref{figura1}, right, we observe that  the absolute error of the approximate solution, $y^h$, is very small, namely, the maximum of the absolute error is approximately $6\times 10^{-11}$, even with a stepsize not too small, $h=1/20$.  This is not surprising,
once the solution belongs to $V_{2}^{1/2}=\left\langle 1,t^{1/2},t,t^{3/2}\right\rangle$.

\begin{table}[h]

\tabcolsep7pt
\begin{adjustbox}{width=1\textwidth}

\begin{tabular}{|l|ccc|ccc|ccc|} \hline
  & \multicolumn{3}{ |c|}{$\alpha=1/4$}&
  \multicolumn{3}{ |c|}{$\alpha=1/2$ } & \multicolumn{3}{ |c|}{$\alpha=2/3$ }
  \\[6pt]
\cline{2-10}
$h$& $\ds y(0)$ &$\ds e^h_{\epsilon}(0.5)$ & $\ds  e^h_{\epsilon}(1)$ &$\ds y(0)$ &$\ds e^h_{\epsilon}(0.5)$ & $\ds  e^h_{\epsilon}(1)$
&$\ds y(0)$ &$\ds e^h_{\epsilon}(0.5)$ & $\ds  e^h_{\epsilon}(1)$
    \\[6pt]
\hline
$\ds 1/10$ & $-0.001564 $& $2.275\times 10^{-11}$&$ 6.382\times 10^{-4}$& $-0.005128 $& $1.975\times 10^{-11}$&$2.297\times 10^{-3}$& $-0.009602 $& $1.279\times 10^{-11}$&$4.604\times 10^{-3}$
  \\
$\ds 1/20$ & $-0.000510 $& $2.938\times 10^{-11}$&$ 2.045\times 10^{-4}$ & $-0.001906 $& $2.844\times 10^{-11}$&$8.438\times 10^{-4}$
& $-0.003922 $& $2.719\times 10^{-11}$&$1.866\times 10^{-3}$

\\
$\ds 1/40$ & $-0.000163 $& $9.517\times 10^{-12}$&$ 6.446\times 10^{-5}$& $-0.000698 $& $1.745\times 10^{-11}$&$3.063\times 10^{-4}$ & $-0.001587 $& $2.411\times 10^{-11}$&$7.510\times 10^{-4}$

\\
$\ds 1/80$ & $ -0.000051 $& $3.158\times 10^{-11}$&$ 2.006\times 10^{-5}$ & $-0.000253 $& $3.025\times 10^{-11}$&$1.103\times 10^{-4}$

& $-0.000638 $& $3.731\times 10^{-11}$&$3.007\times 10^{-4}$
\\
$\ds 1/160$ & $-0.000016 $& $4.922\times 10^{-12}$&$ 6.189\times 10^{-6}$ & $ -0.000032
 $& $ 8.683\times 10^{-11}$&$3.948\times 10^{-4}$

& $-0.000255 $& $2.321\times 10^{-11}$&$1.200\times 10^{-4}$
 \\
$\ds 1/320$ &$ -0.000005 $& $4.251\times 10^{-12}$&$ 1.896\times 10^{-6}$ & $-0.000091 $& $1.230\times 10^{-11}$&$1.408\times 10^{-4}$
& $-0.000102 $& $1.198\times 10^{-11}$&$4.781\times 10^{-5}$
 \\
\hline
\end{tabular}
 \end{adjustbox}
 \caption{Example \ref{exemplo1} with several values of $\alpha$. Comparison with the exact solution at $t=0.5$ (the value that defines the boundary condition) and $t=1$ with several values of the stepsize $h$.} \label{tab1}
\end{table}


\begin{table}[h]

 \centering
\begin{tabular}{|l|cc|cc|cc|} \hline
  & \multicolumn{2}{ |c|}{$\alpha=1/4$}&
  \multicolumn{2}{ |c|}{$\alpha=1/2$ } & \multicolumn{2}{ |c|}{$\alpha=2/3$ }
  \\[6pt]
\cline{2-7}
$N$& $\ds \|e^{a/N}\|_{\infty}$ &$\ds p_N$ & $\ds \|e^{a/N}\|_{\infty}$ &$\ds p_N$&$\ds \|e^{T/N}\|_{\infty}$ &$\ds p_N$
    \\[6pt]
\hline

$\ds 10$ & $1.564\times 10^{-3} $& $ -$& $5.128\times 10^{-3} $& $- $& $9.602\times 10^{-3} $& $- $
  \\
$\ds 20$& $5.098\times 10^{-4} $& $1.62 $& $1.906\times 10^{-3} $& $1.43 $& $3.922\times 10^{-3} $& $1.29 $
\\
$\ds 40$& $1.626\times 10^{-4} $& $ 1.65$& $6.978\times 10^{-4} $& $1.45 $& $1.587\times 10^{-3} $& $1.31 $
\\
$\ds 80$& $5.107\times 10^{-5} $& $1.67 $& $2.527\times 10^{-4} $& $1.47 $& $6.381\times 10^{-4} $& $1.31 $
\\
$\ds 160$& $1.586\times 10^{-5} $& $1.69 $& $9.084\times 10^{-5} $& $1.48 $& $2.553\times 10^{-4} $& $1.32 $
\\
$\ds 320$& $4.883\times 10^{-6} $& $1.70 $& $3.249\times 10^{-5} $& $1.48 $& $1.019\times 10^{-4} $& $1.33 $
\\
\hline

\end{tabular}

 \caption{Example \ref{exemplo1} with several values of $\alpha$. Maximum of absolute errors and experimental order of convergence.}\label{tab2}
\end{table}

\begin{table}[h]
\tabcolsep7pt
\begin{adjustbox}{width=1\textwidth}

\begin{tabular}{|l|ccc|ccc|ccc|} \hline
  & \multicolumn{3}{ |c|}{$\alpha=1/4$}&
  \multicolumn{3}{ |c|}{$\alpha=1/2$ } & \multicolumn{3}{ |c|}{$\alpha=2/3$ }
  \\[6pt]
\cline{2-10}
$h$& $\ds y(0)$ &$\ds e^h_{\epsilon}(0.5)$ & $\ds  e^h_{\epsilon}(1)$ &$\ds y(0)$ &$\ds e^h_{\epsilon}(0.5)$ & $\ds  e^h_{\epsilon}(1)$
&$\ds y(0)$ &$\ds e^h_{\epsilon}(0.5)$ & $\ds  e^h_{\epsilon}(1)$
    \\[6pt]
    \hline
$\ds 1/20$ & $9.313 \times 10^{-11} $&$3.014  \times 10^{-3}$& $5.280 \times 10^{-3}$& $9.313 \times 10^{-11}  $& $8.484\times 10^{-4}$&$1.487\times 10^{-3}$& $9.313 \times 10^{-11}  $& $3.440  \times 10^{-4}$&$6.534  \times 10^{-4}$
  \\
$\ds 1/40$  & $9.313 \times 10^{-11} $& $1.364  \times 10^{-3}$ & $2.386 \times 10^{-3}$ & $9.313 \times 10^{-11}  $& $3.273\times 10^{-4}$&$5.559\times 10^{-4}$& $9.313 \times 10^{-11}  $& $1.214  \times 10^{-4}$ &$  2.181\times 10^{-4}$
\\
$\ds 1/80$ &   $9.313 \times 10^{-11} $&$5.908  \times 10^{-4}$& $1.038 \times 10^{-3}$ & $9.313 \times 10^{-11}  $& $1.214\times 10^{-4}$&$2.028\times 10^{-4}$& $9.313 \times 10^{-11}  $& $4.121  \times 10^{-5}$&$  7.127\times 10^{-5}$
\\
$\ds 1/160$ & $9.313 \times 10^{-11} $&$2.504   \times 10^{-4}$& $4.427 \times 10^{-4}$ & $9.313 \times 10^{-11}  $& $4.412\times 10^{-5}$&$7.302\times 10^{-5}$& $9.313 \times 10^{-11}  $& $1.369  \times 10^{-5}$ &$  2.302\times 10^{-5}$

\\
$\ds 1/320$  & $9.313 \times 10^{-11} $&$1.051 \times 10^{-4}$& $1.872 \times 10^{-4}$ & $9.313 \times 10^{-11}  $& $1.587\times 10^{-5}$&$
2.611\times 10^{-5}$ & $9.313 \times 10^{-11}  $ & $4.486  \times 10^{-6}$ &  $  7.382\times 10^{-6}$
 \\
\hline
\end{tabular}
 \end{adjustbox}
 \caption{Example \ref{exemplo2} with several values of $\alpha$. Comparison with the exact solution at $t=0.5$ (the value in  the boundary condition) and $t=1$ with several values of the stepsize $h$ (shooting method with Method 2 to
solve the IVP).}\label{tab3}
\end{table}


\begin{table}[h]
\centering
\begin{tabular}{|l|cc|cc|cc|} \hline
  & \multicolumn{2}{ |c|}{$\alpha=1/4$}&
  \multicolumn{2}{ |c|}{$\alpha=1/2$ } & \multicolumn{2}{ |c|}{$\alpha=2/3$ }
  \\[6pt]
\cline{2-7}
$N$& $\ds \|e^{T/N}\|_{\infty}$ &$\ds p_N$ & $\ds \|e^{T/N}\|_{\infty}$ &$\ds p_N$&$\ds \|e^{T/N}\|_{\infty}$ &$\ds p_N$
    \\[6pt]
\hline

$\ds 20$ & $5.306\times 10^{-3} $& $-$ & $1.494\times 10^{-3} $& $ -$ & $6.565\times 10^{-4} $& $- $
  \\
$\ds 40$& $2.399\times 10^{-3} $& $1.15 $ & $5.592\times 10^{-4} $& $ 1.42$& $2.194\times 10^{-4} $& $1.58$
\\
$\ds 80$& $1.043\times 10^{-3} $& $1.20 $& $2.041\times 10^{-4} $& $1.45$& $7.182\times 10^{-5} $& $1.61$
\\
$\ds 160$& $4.450\times 10^{-4} $& $1.23 $& $7.352\times 10^{-5} $& $ 1.47$& $2.322\times 10^{-5} $& $1.63$
\\
$\ds 320$& $1.881\times 10^{-4} $& $1.24$ & $2.630\times 10^{-5} $& $1.48$& $7.452\times 10^{-6} $& $1.64 $
\\
\hline

\end{tabular}

 \caption{Example \ref{exemplo2} with several values of $\alpha$. (shooting method with Method 2 to
solve the IVP).}\label{tab4}
\end{table}

In Figures \ref{figura1} the absolute  errors of the approximate solutions of Examples \ref{exemplo1}  and \ref{exemplo2}, with several values of $\alpha$, are plotted for  stepsize $\ds h=1/160$. For example  \ref{exemplo1} we observe that the absolute error is minimum at the point $t=a$ and for example  \ref{exemplo2}  the absolute error is minimum at the point $t=0$. For both examples the absolute error decreases with the value of $\alpha$.

\begin{figure}[h]
\centering
\begin{tabular}{cc}
\includegraphics[scale=0.5]
{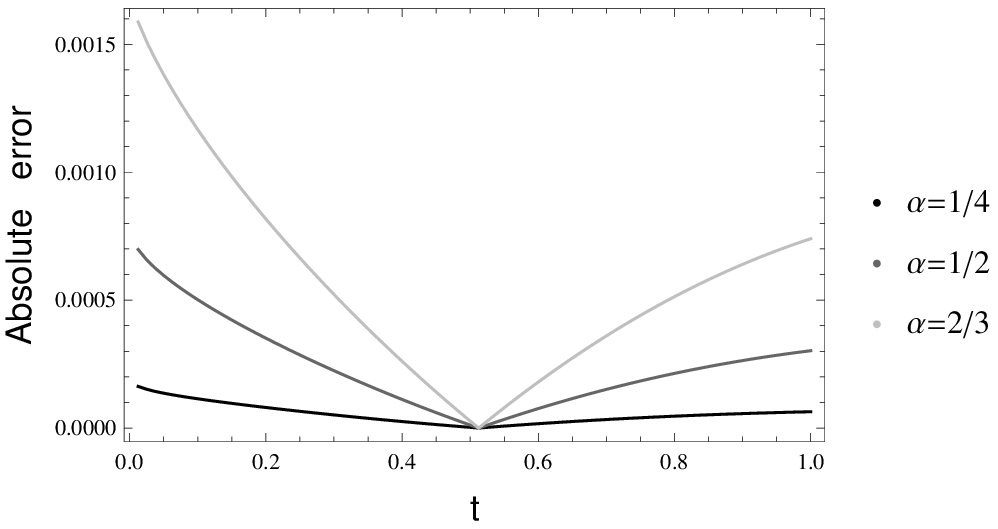}&\includegraphics[scale=0.5]
{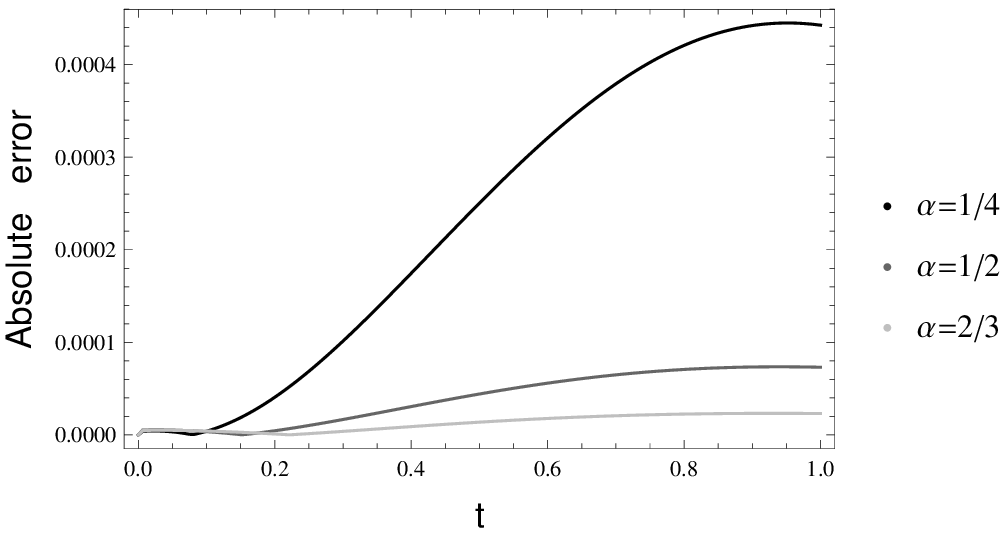}
\end{tabular}
\caption{
 Plot of error function $|y(t)-y^{h}(t)|$, with $h=1/160$, for Example
\ref{exemplo1} (left)  for Example
\ref{exemplo2} (right), with $\alpha=1/4$, $\alpha=1/2$ and $\alpha=2/3$. }
\label{figura1}
\end{figure}

\begin{table}[h]

\tabcolsep7pt
\begin{adjustbox}{width=1\textwidth}

\begin{tabular}{|l|ccc|ccc|} \hline
  & \multicolumn{3}{ |c|}{Method 1}&
  \multicolumn{3}{ |c|}{Method 3}
  \\[6pt]
\cline{2-7}
$h$& $\ds y(0)$ &$\ds e^h_{\epsilon}(0.5)$ & $\ds  e^h_{\epsilon}(1)$ &$\ds y(0)$ &$\ds e^h_{\epsilon}(0.5)$ & $\ds  e^h_{\epsilon}(1)$ \\
    \hline
$\ds 1/20$ & $3.037 \times 10^{-3} $&$1.340 \times 10^{-12}$& $1.814 \times 10^{-5}$& $   5.821\times 10^{-11}  $& $7.844\times 10^{-5}$&$6.302\times 10^{-6}$
  \\
$\ds 1/40$ & $1.121 \times 10^{-3} $&$8.904 \times 10^{-12}$& $4.458 \times 10^{-6}$& $ 5.821\times 10^{-11}  $& $2.955\times 10^{-5}$&$1.657\times 10^{-6}$
  \\
$\ds 1/80$& $4.081 \times 10^{-4} $&$1.305 \times 10^{-11}$& $1.105 \times 10^{-6}$& $ 5.821\times 10^{-11}  $& $1.098\times 10^{-5}$&$   4.286\times 10^{-7}$
  \\
$\ds 1/160$ & $1.472 \times 10^{-4} $&$1.531 \times 10^{-11}$& $2.750 \times 10^{-7}$& $ 5.821\times 10^{-11}  $& $ 4.119\times 10^{-7}$&$1.096\times 10^{-7}$
  \\
$\ds 1/320$  & $5.275 \times 10^{-5} $&$6.367 \times 10^{-12}$& $6.859 \times 10^{-8}$& $ 5.821\times 10^{-11}  $& $1.054\times 10^{-7}$&$  2.785\times 10^{-8}$
  \\
\hline
\end{tabular}
 \end{adjustbox}
 \caption{Example \ref{exemplo3} with $\alpha=1/2$. Comparison with the exact solution at $t=0$, $t=0.5$ (the value that defines the boundary condition) and $t=1$ with several values of the stepsize $h$ (shooting method with Method 1 and Method 3 on the space $V_{h\,1}^{1/2}$ ($c_1=1/3$, $c_2=1$) to
solve the IVP).}\label{tab5exemp3}
\end{table}

\bc
\begin{table}[h]
\centering
\begin{tabular}{|l|cc|cc|} \hline
  & \multicolumn{2}{ |c|}{Method 1}&
  \multicolumn{2}{ |c|}{Method 3}
  \\[6pt]
\cline{2-5}
$N$& $\ds \|e^{T/N}\|_{\infty}$ &$\ds p_N$ & $\ds \|e^{T/N}\|_{\infty}$ &$\ds p_N$
    \\[6pt]
\hline
$\ds 20$ & $3.304\times 10^{-3} $& $ -$& $7.844\times 10^{-5} $& $- $
  \\
$\ds 40$ & $1.121\times 10^{-3} $& $ 1.44$& $2.955\times 10^{-5} $& $1.41 $
  \\
$\ds 80$ & $4.081\times 10^{-4} $& $ 1.46$& $1.098\times 10^{-5} $& $1.43 $
  \\
$\ds 160$ & $1.472\times 10^{-4} $& $ 1.47$& $3.981\times 10^{-6} $& $1.46 $
  \\
$\ds 320$ & $5.275\times 10^{-5} $& $ 1.48$& $1.425\times 10^{-6} $& $1.48$
  \\
\hline

\end{tabular}
 \caption{
Example \ref{exemplo3} with $\alpha=1/2$.
  Maximum of absolute errors and experimental order of convergence  (shooting method with Method 1 and Method 3 on the space $V_{h\,1}^{1/2}$ ($c_1=1/3$, $c_2=1$).}
 \label{tab6exemp3}
\end{table}
\ec

\begin{figure}[h]
\centering
\begin{tabular}{cc}
\includegraphics[scale=0.5]
{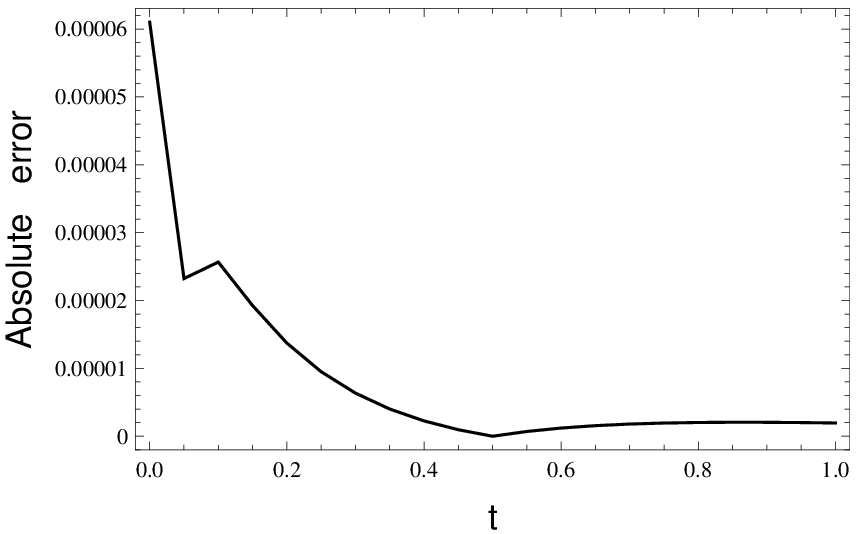}&\includegraphics[scale=0.5]
{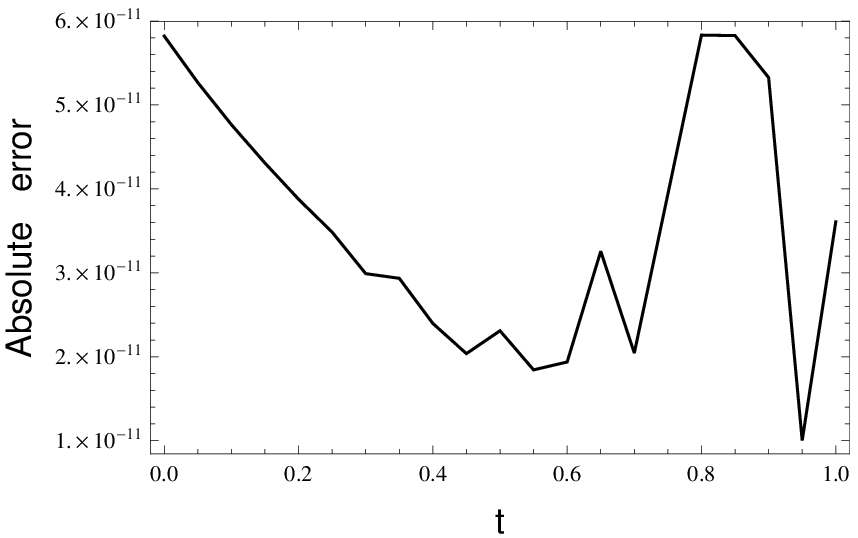}
\end{tabular}
\caption{
 Plot of error function $|y(t)-y^{h}(t)|$, with $h=1/20$, from Example
\ref{exemplo3}. Left: Shooting method with Method 3 on the space $V_{h\,1}^{1/2}$.
Right: Shooting method with Method 3 on the space $V_{h\,2}^{1/2}$.  }
\label{figura1}
\end{figure}

\clearpage

 \subsection{Dependence on the problem parameters }
 \noindent

In this subsection we consider a nonlinear problem and illustrate numerically the stability of the problem.

Let us consider the tempered fractional differential equation
\begin{ex}\label{exemplo4}
\be\nonumber && \mathbb{D}^{\alpha,\lambda}
\left(y(t)\right)=2t+\frac{\Gamma(\alpha+1)}{3\exp(\lambda a)a^{\alpha}}\sin(u)=f(t,u),\quad t>0,\\
\nonumber &&
y(a)=y_a,
\ee
\end{ex}
with $\lambda=2$, $a=1/2$, $y_a=1$ and $\alpha=1/2$. Note that the function $f$ satisfies the assumptions of Theorems \ref{Varua}-\ref{Varf}. In this case the exact solution of Example \ref{exemplo4} is unknown.
\\ Let us consider the perturbed problems

\be &&
\label{expertb1}
 \mathbb{D}^{\alpha,\lambda}
\left(z(t)\right)=f(t,z),\quad t>0,\\
\nonumber &&
z(a)=y_a+\epsilon_{bc},
\\
\nonumber &&\\
&&\label{expertb2}
 \mathbb{D}^{\alpha,\lambda}
\left(z(t)\right)=f(t,z)+\epsilon_{f},\quad t>0,\\
\nonumber &&
z(a)=y_a,
\\
\nonumber &&\\
&&\label{expertb3}
 \mathbb{D}^{\alpha,\lambda+\epsilon_{\lambda}}
\left(z(t)\right)=f(t,z),\quad t>0,\\
\nonumber && z(a)=y_a,\\
\nonumber && \\
&&\label{expertb4}
 \mathbb{D}^{\alpha+\epsilon_{\alpha},\lambda}
\left(z(t)\right)=f(t,z),\quad t>0,\\
\nonumber &&
z(a)=y_a.
\ee

The obtained $\ds\max_{1\leq i \leq N} \|y_i-z_i\|=\|y-z\|_{\infty}$ are presented in Tables \ref{tab7}, \ref{tab8} and \ref{tab9}, where
$y_i$ and $z_i$ are the obtained numerical approximations of $y(t)$ and $z(t)$ at the discretization
points $t=t_i=ih$, with $h=a/N$, and $z$ is the solution of the perturbed problems (\ref{expertb1}), (\ref{expertb2}) and (\ref{expertb3}), respectively.\\   In Table \ref{tab7} we present the results obtained when we compare the problems (\ref{exemplo4}) and (\ref{expertb1}), when the  boundary condition suffers a perturbation.

In Table \ref{tab8} we present the results obtained when we compare the problems (\ref{exemplo4}) and (\ref{expertb2}), when the  source function $f$ has a perturbation, $\epsilon_f$.

Finally, in Table \ref{tab9} we illustrate how the solution of (\ref{expertb3}) varies  with $\ds \epsilon_{\lambda}$.

 According to the numerical results in Tables \ref{tab7}, \ref{tab8}  and \ref{tab9}, we see that, independently of the used step
size $h$, we have  $\ds \|y-z\|_{\infty} \sim \epsilon_{bc}$,
$\ds \|y-z\|_{\infty} \sim \epsilon_{f}$ and $\ds \|y-z\|_{\infty} \sim \epsilon_{\lambda}$, if $z$ is the approximate solution of the problems (\ref{expertb1}), (\ref{expertb2}) and (\ref{expertb3}), respectively. The numerical results are in agreement with the theoretical results proved in Theorems \ref{Varua}, \ref{Varlambda} and \ref{Varf}.

\begin{table}[h]

\tabcolsep7pt
\begin{adjustbox}{width=1\textwidth}

\begin{tabular}{|l|ccccc|}
\hline
  & \multicolumn{5}{ |c|}{Values of $\ds \epsilon_{bc}$}  \\
\cline{2-6}
$h$& $\ds 0.1$ &$\ds 0.01$& $\ds 0.001 $&$\ds 0.0001 $&$\ds 0.00001 $
  \\
\hline
$\ds 1/20$ & $2.5704\times 10^{-1 }$& $2.5518\times 10^{-2 }$& $2.5499\times 10^{-3 }$& $2.5498\times 10^{-4 }$& $2.5497\times 10^{-5 }$\\
$\ds 1/40$ & $2.5712\times 10^{-1 }$& $2.5525\times 10^{-2 }$& $2.5507\times 10^{-3 }$& $2.5506\times 10^{-4 }$& $2.5505\times 10^{-5 }$\\
$\ds 1/80$ & $2.5715\times 10^{-1 }$& $2.5528\times 10^{-2 }$& $2.5510\times 10^{-3 }$& $2.5508\times 10^{-4 }$& $2.5508\times 10^{-5 }$\\
$\ds 1/160$ & $2.5716\times 10^{-1 }$& $2.5529\times 10^{-2 }$& $2.5519\times 10^{-3 }$& $2.5509\times 10^{-4 }$& $2.5509\times 10^{-5 }$
\\
\hline
\end{tabular}
 \end{adjustbox}
 \caption{ Maximum of the asbolute errors,
 $|y^h-z^h|$, where $y^h$ is  the  numerical solution of problem (\ref{exemplo4}) and $z^h$ the numerical solution of the problem  (\ref{expertb1}) with several values of $\epsilon_{bc}$.}\label{tab7}
\end{table}

\begin{table}[h]

\tabcolsep7pt
\begin{adjustbox}{width=1\textwidth}

\begin{tabular}{|l|ccccc|}
\hline
  & \multicolumn{5}{ |c|}{Values of $\ds \epsilon_{f}$}  \\
\cline{2-6}
$h$& $\ds 0.1$ &$\ds 0.01$& $\ds 0.001 $&$\ds 0.0001 $&$\ds 0.00001 $
  \\
\hline
$\ds 1/20$ & $7.8815 \times 10^{-2 }$& $7.8841\times 10^{-3 }$& $7.8844 \times 10^{-4 }$& $7.8844 \times 10^{-5 }$& $  7.8845\times 10^{-6 }$\\
$\ds 1/40$ & $7.8818 \times 10^{-2 }$& $7.8844\times 10^{-3 }$& $7.8847 \times 10^{-4 }$& $7.8847 \times 10^{-5 }$& $  7.8847\times 10^{-6 }$\\
$\ds 1/80$& $7.8819 \times 10^{-2 }$& $7.8845\times 10^{-3 }$& $7.8848 \times 10^{-4 }$& $7.8848 \times 10^{-5 }$& $  7.8848\times 10^{-6 }$\\
$\ds 1/160$ & $7.8820 \times 10^{-2 }$& $7.8845\times 10^{-3 }$& $7.8848 \times 10^{-4 }$& $7.8848 \times 10^{-5 }$& $  7.8848\times 10^{-6 }$\\
\hline
\end{tabular}
 \end{adjustbox}
 \caption{ Maximum of the absolute errors,
 $|y^h-z^h|$, where $y^h$ is  the  numerical solution of problem (\ref{exemplo4}) and $z^h$ the numerical solution of the problem  (\ref{expertb2}) with several values of $\epsilon_{f}$.
 }\label{tab8}
\end{table}

\begin{table}[h]

\tabcolsep7pt
\begin{adjustbox}{width=1\textwidth}

\begin{tabular}{|l|ccccc|}
\hline
  & \multicolumn{5}{ |c|}{Values of $\ds \epsilon_{\lambda}$}  \\
\cline{2-6}
$h$& $\ds 0.1$ &$\ds 0.01$& $\ds 0.001 $&$\ds 0.0001 $&$\ds 0.00001 $
  \\
\hline
$\ds 1/20$ & $2.2715 \times 10^{-1 }$& $2.1483\times 10^{-2 }$& $2.1364 \times 10^{-3 }$& $2.1352 \times 10^{-4 }$& $  2.1351\times 10^{-5 }$\\
$\ds 1/40$ & $2.2726 \times 10^{-1}$& $2.1493\times 10^{-2 }$& $2.1374 \times 10^{-3 }$& $2.1362 \times 10^{-4 }$& $  2.1361\times 10^{-5 }$\\
$\ds 1/80$& $2.2729 \times 10^{-1 }$& $2.1496\times 10^{-2 }$& $2.1377 \times 10^{-3 }$& $ 2.1365\times 10^{-4 }$& $  2.1364\times 10^{-5 }$\\
$\ds 1/160$ & $2.2730 \times 10^{-1 }$& $2.1497\times 10^{-2 }$& $2.1378 \times 10^{-3 }$& $ 2.1366\times 10^{-4 }$& $  2.1365\times 10^{-5 }$
\\
\hline
\end{tabular}
 \end{adjustbox}
 \caption{ Maximum of the asbolute errors,
 $|y^h-z^h|$, where $y^h$ is  the  numerical solution of problem (\ref{exemplo4}) and $z^h$ the numerical solution of the problem  (\ref{expertb3}) with several values of $\epsilon_{\lambda}$.
}\label{tab9}
\end{table}

\begin{figure}[h]
\centering
\begin{tabular}{ccc}
\includegraphics[scale=0.35]
{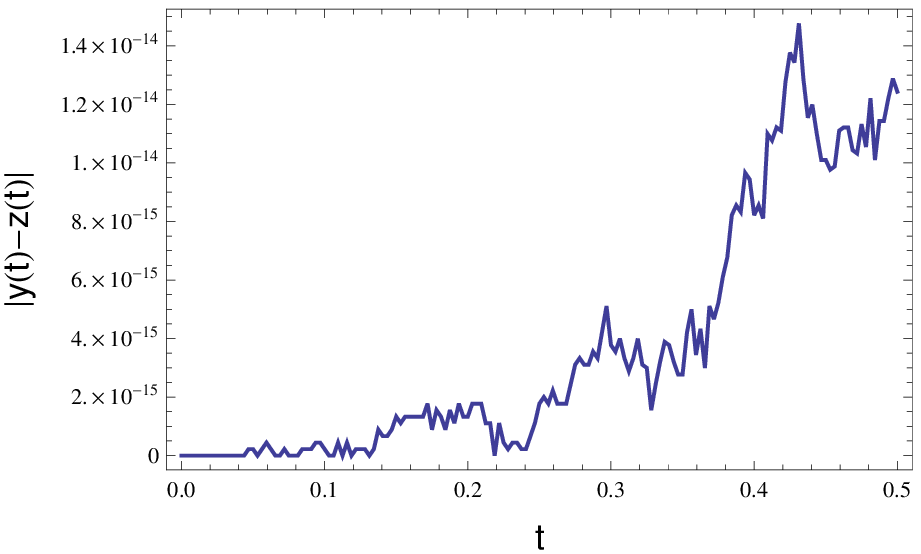}&
\includegraphics[scale=0.34]
{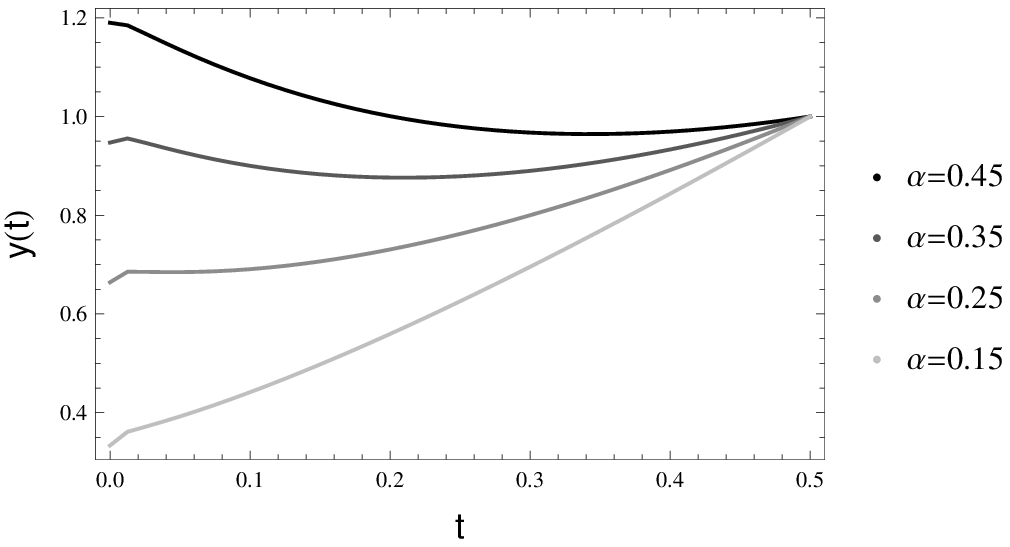} &\includegraphics[scale=0.34]
{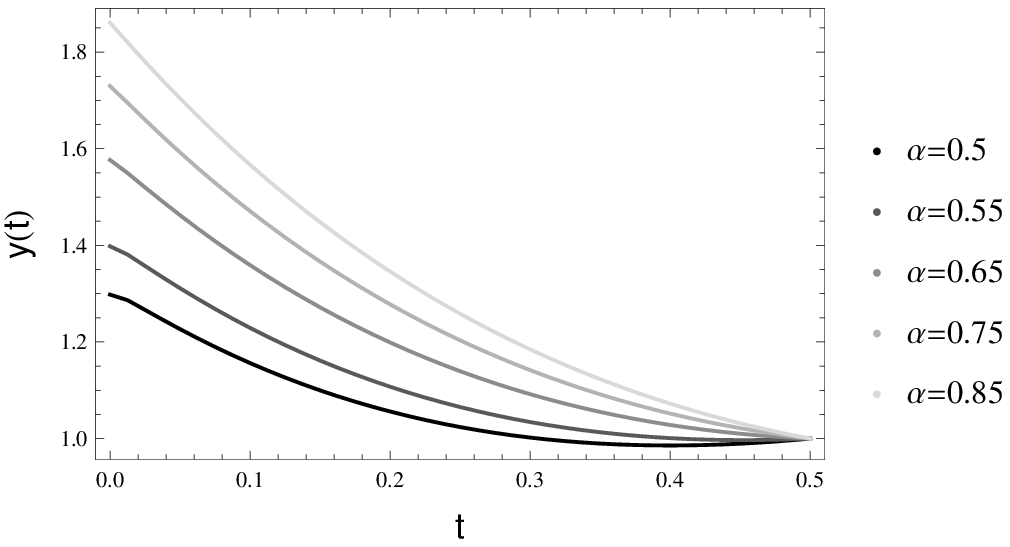}
\end{tabular}
\caption{
Left: Plot of the error function $|y(t)-z(t)|$, where $z$ is the approximate solution of
(\ref{expertb4})  with $\epsilon_{\alpha}=0.01$. The approximate solutions are obtained using the Method 2  with $\ds h=1/160$.
Center and right: Approximate solutions of (\ref{expertb4}) with several values of $\alpha$. }
\label{figuraVaralpha}

\end{figure}

\begin{figure}[h]
\centering
\begin{tabular}{cc}
\includegraphics[scale=0.4]
{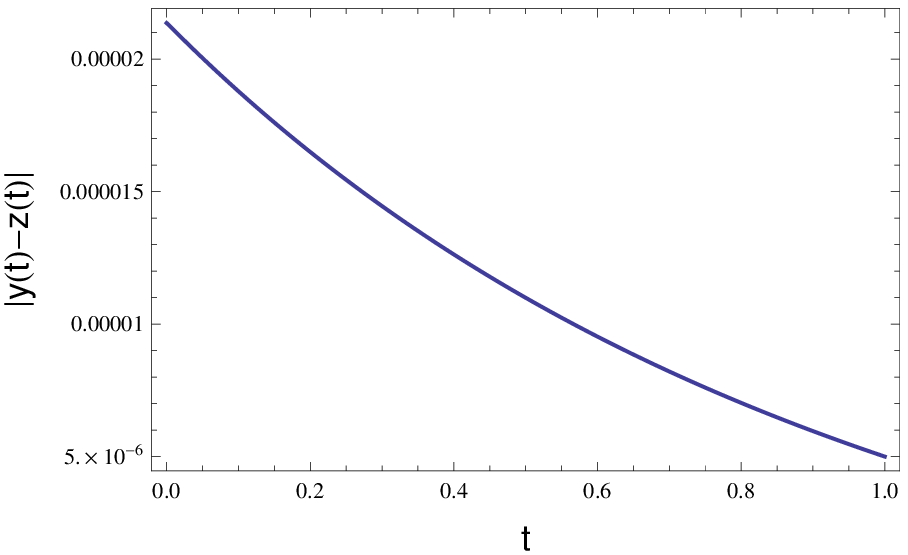} &\includegraphics[scale=0.36]
{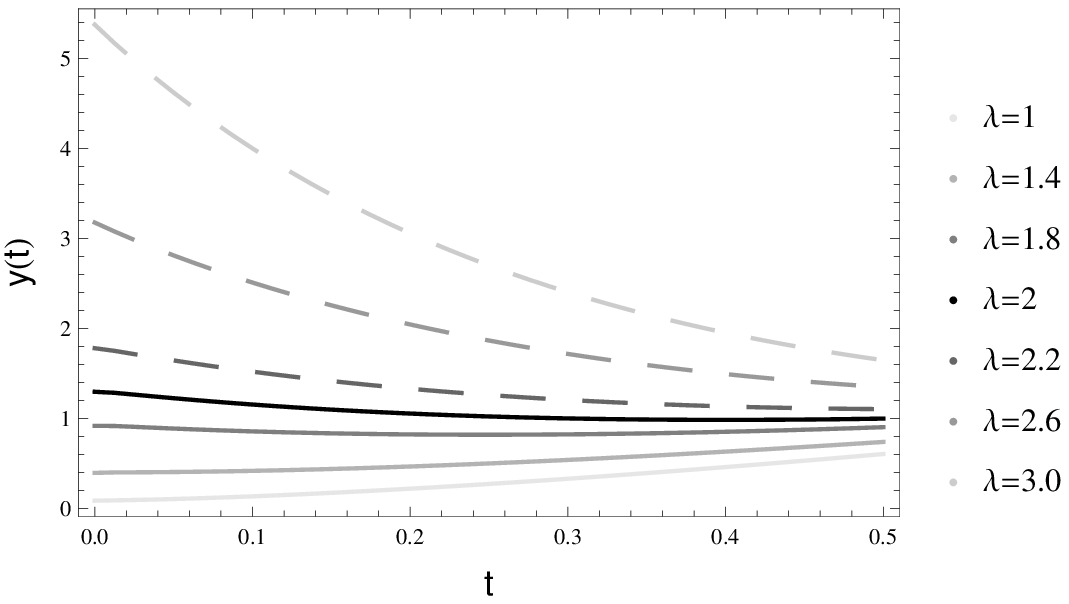}
\end{tabular}
\caption{
Left: Plot of the error function $|y(t)-z(t)|$, where $z$ is the approximate solution of the (\ref{expertb3})  with $\epsilon_{\alpha}=0.00001$. The approximate solutions are obtained using the Method 2  with $\ds h=1/160$.
Right: Approximate solutions of (\ref{expertb3}) with several values of
$\lambda$. The plots with dashed lines are related with the values of $\lambda$ greater than 2.  }
\label{figuraVarLambda}
\end{figure}

In Figure \ref{figuraVaralpha} we present an approximate solution of the problem (\ref{expertb3}) with $\epsilon_\alpha=0.01$, and we observe that the variation is very small. We also  plot the approximate solution of  (\ref{exemplo4}), for several values of $\alpha$, and we observe that the solution is an increasing function for $\alpha<0.5$ and a decreasing function for $\alpha\geq 0.5$.   Finally, in Figure \ref{figuraVarLambda} we plot the absolute error
 $|y^{1/160}-z^{1/160}|$, where $y^{1/160}$ is  the  approximate  solution of problem (\ref{exemplo4}) and $z^{1/160}$ the approximate  solution of the problem  (\ref{expertb3}) with $\epsilon_{\lambda}=10^{-5}$. It can be observed that the absolute error is less than $\lambda \times 10^{-5}$ and the absolute error is maximum at the origin.

\section{Conclusions}
\
\noindent

We have analysed the well-posedness of ordinary tempered terminal value problems. Based on the relationship between non-tempered and tempered Caputo derivatives we have proposed three numerical schemes to approximate the solution of such problems. It should be noted that Method 3 has the advantage to properly deal with nonsmooth solutions which constitutes an important feature in the numerical approximation of fractional differential problems. In the future, we intend to extend it to partial and distributed differential problems.

\vspace{0.3cm}

\section* {Acknowledgments}
\
\noindent

 The two authors acknowledge financial support from FCT – Funda\c c\~{a}o para a Ci\^{e}ncia e a Tecnologia (Portuguese Foundation for Science and Technology), through Project UID/MAT/00013/2013 and project UID/MAT/00297/2013, respectively.

\end{document}